\documentclass[12pt]{article}
\usepackage{amssymb,amsmath,graphicx}
\usepackage{longtable}
\usepackage{hyperref}
                               
\newcommand{\tab}{\hspace{5mm}}

\begin{document}

\begin{center}
\textbf{{\huge G\"{O}DEL'S THEOREM IS\\
\vspace{0.2 cm}
INVALID}}\\
\vspace{0.5 cm}
{\Large Diego Sa\'{a}}
\footnote{Escuela Polit\'{e}cnica Nacional. Quito --Ecuador. email: dsaa@server.epn.edu.ec}\\
Copyright {\copyright}2000\\
\end{center}

\begin{abstract}
{G\"{o}del's results have had a great impact in diverse fields such 
as philosophy, computer sciences and fundamentals of mathematics.\\

The fact that the rule of mathematical induction is contradictory 
with the rest of clauses used by G\"{o}del to prove his undecidability 
and incompleteness theorems is proved in this paper. This means that those theorems are invalid.\\

In section 1, a study is carried out on the mathematical induction {\underline {principle}}, 
even though it is not directly relevant to the problem, just 
to familiarize the reader with the operations that are used later; 
in section 2 the {\underline {rule}} of mathematical induction is introduced, 
this rule has a metamathematical character; in section 3 the 
original proof of G\"{o}del's undecidability theorem is reproduced, 
and finally in section 4 the same proof is given, but now with 
the explicit and formal use of all the axioms; this is needed 
to be able to use logical resolution. It is shown that the inclusion 
of the mathematical induction rule causes a contradiction.}
\end{abstract}


\section{INTRODUCTION}

\indent Even though mathematics itself has never lost its power, the 
limitative theorems need the conceptual admission that there 
are unanswerable questions. As Howard Delong \cite{delong} relates, ``previously 
it was thought that if a question was well-defined, that question 
had an answer''. An illustration of this attitude may be 
found in Hilbert's address ``on the infinite'', delivered 
in 1925: ``As an example of the way in which fundamental 
questions can be treated I would like to choose the thesis that 
every mathematical problem can be solved. We are all convinced 
of that. After all, one of the things that attract us most when 
we apply ourselves to a mathematical problem is precisely that 
within us we always hear the call: here is the problem, search 
for the solution; you can find it by pure thought, for in mathematics 
there is no ignorabimus.''\\

In the same way, Nagel and Newman \cite{nagel} in his book published in 
1958 confirm that ``until recently it was taken as a matter 
of course that a complete set of axioms for any branch of mathematics 
can be assembled. In particular, mathematicians believed that 
the set proposed for arithmetic in the past was in fact complete, 
or, at worst, could be made complete simply by adding a finite 
number of axioms to the original list. The discovery that this 
will not work is one of G\"{o}del's major achievements''.\\

In another place of their book, Nagel and Newman \cite{nagel} evaluate 
the effects of those theorems declaring that: ``the import 
of G\"{o}del's conclusions is far reaching, though it has not been 
fully fathomed. These conclusions show that the prospect of finding 
for every deductive system (and in particular, for a system in 
which the whole arithmetic can be expressed) an absolute proof 
of consistency that satisfies the finitistic requirements of 
Hilbert's proposal, though not logically impossible, is most 
unlikely''.\\

There are consequences in several fields, from philosophy to 
computer sciences, but the above mentioned are enough to justify 
the need to make sure that G\"{o}del's theorems are free of any 
kind of faults.\\

In the present paper a contradiction is exposed, between those 
theorems and an accepted and well known rule of mathematics, 
which allows us to conclude that such incompleteness and undecidability 
theorems are invalid.\\

In section 1, a revision is made of the mathematical induction {\underline {principle}}, 
even though it is not directly relevant to the problem, just 
to get the reader acquainted with the operations used later; 
the {\underline {rule}} of mathematical induction, which has a metamathematical 
character, is introduced in section 2; in section 3 the proof 
of G\"{o}del's undecidability theorem is exhibited again but now 
with the explicit use of all the axioms, which are used with 
the objective of employing logical resolution; the further inclusion 
of the mathematical induction rule causes the cited contradiction.

\section{THE MATHEMATICAL INDUCTION\\
PRINCIPLE}

The mathematical induction principle can be expressed in the 
following way:\\

Let P be any property of natural numbers.\\

Suppose that:\\

A) 0 (zero) has the property P.\\

B) If any natural number x has property P, then its successor 
s(x) also has property P.\\

Then:\\

C) All natural numbers x have property P.\\

This ``mathematical induction principle'' can be 
expressed more concisely by using the logical connectives for 
conjunction (\ensuremath{\wedge} ) and implication (\ensuremath{\rightarrow} ) and 
the universal quantifier (\ensuremath{\forall} ), as follows:\\

\begin{equation}
P(0) \wedge (\forall x)( P(x) \rightarrow P(s(x)) ) \rightarrow 
(\forall x) P(x) 
\end{equation}

This formula is referred to by some authors as an \textit{axiom}, 
because it corresponds to the fifth Peano's axiom for the numbers 
(See G\"{o}del \cite{goedel}, p.61). In the APPENDIX of this paper the author presents a prove that it is a \textit{theorem}.\\

Let us work for a moment with premises (A) and (B) only:

\begin{equation}
P(0)
\label{eq2}
\end{equation}

\begin{equation}
(\forall x)( P(x) \rightarrow P(s(x)) )
\label{eq3}
\end{equation}\\

If we assign to x the value 0 in (\ref{eq3}), it can be proved, by means 
of a truth table, that (\ref{eq2}) and (\ref{eq3}) are logically equivalent to:

\begin{equation}
P(0)\wedge P(s(0))
\end{equation}

With a new value for x in (\ref{eq3}), now of s(0), clauses (\ref{eq2}) and (\ref{eq3}) 
are equivalent to:

\begin{equation}
P(0)\wedge P(s(0)) \wedge P(s(s(0)))
\label{eq5}
\end{equation}

Expression (\ref{eq5}) is an expansion of:

\begin{equation}
(\forall x) P(x)
\end{equation}

where the universal quantifier ``for all x'' refers 
to the elements of the set: 

\begin{equation}
\{\ 0, s(0), s(s(0))\ \}
\end{equation}

The same operations could be successively applied to all the 
integers, with the result that the premises (A) and (B) are equivalent 
to $(\forall x) P(x)$, this is identical to the conclusion 
(C). So the mathematical induction principle is logically equivalent to the 
clause:

\begin{equation}
(\forall x) P(x) \rightarrow (\forall x) P(x)
\end{equation}

which is a tautology, and as such it does not provide any new 
useful information to number theory, so it is superfluous and 
can be removed.\\

\section{THE MATHEMATICAL INDUCTION\\
RULE}

The previous conclusion does not mean, however, that Peano didn't 
have something to say. But, in order to express that ``something'', 
it is not enough a mathematical principle but a metamathematical 
rule. In fact, Kleene \cite{kleene} has the ``induction rule'' 
expressed as:

\begin{equation}
\Gamma \vdash P(0),\ and\ \Gamma 
,P(x) \vdash P(s(x)),\ then\ \Gamma \vdash P(x)
\end{equation}

Where the turnstile symbol ``$\vdash$'' represents the metamathematical concept ``then it is deducible''.\\

Let us use G\"{o}del's predicate ``Bew'', which means 
``deducible formula'', instead of $\vdash$. Kleene's 
induction rule can then be reexpressed as:\\

$Bew(P(0)) \wedge (\forall x) ( Bew(P(x)) \rightarrow 
Bew(P(s(x))) )$
\begin{equation} 
\rightarrow Bew( (\forall x) P(s(x)) ) 
\label{eq10}
\end{equation}

If the value 0 is assigned to x, as we did in the case of the induction axiom, it can be found that the antecedent of the implication (\ref{eq10}):\\

$Bew(P(0)) \wedge ( Bew(P(0)) \rightarrow Bew(P(s(0))) )$\\

is equivalent to $Bew(P(0)) \wedge Bew(P(s(0)))$, and in this 
fashion, continuing successively with these operations, it is 
found that the antecedent rule (\ref{eq10}) is logically equivalent to:

\begin{equation}
(\forall x) Bew(P(x))\\
\end{equation}

After each operation, the universal quantifier has a range in 
a set that contains one additional successor than the previous 
set, and eventually could be applied to a set that has whatever 
number of elements.\\

Therefore, Kleene's induction rule, (\ref{eq10}), is logically equivalent 
to:

\begin{equation}
(\forall x) Bew(P(x)) \rightarrow Bew( (\ensuremath{\forall} x) P(x) )
\label{eq12}
\end{equation}

Where the quantifiers on each side of the implication should 
have their ranges within the same set of elements.\\
The implication in the opposite sense:

\begin{equation}
Bew( (\forall x) P(x) ) \rightarrow (\forall x) Bew(P(x))
\label{eq13}
\end{equation}

was used by G\"{o}del, implicitly, in the first part of his proof 
of the undecidability theorem (\cite{goedel}, p.76). It also has its non-metamathematical 
counterpart among the universally valid rules of formal systems 
with the name of ``specialization rule'':

\begin{equation}
(\forall x) P(x) \rightarrow P(c)
\label{eq14}
\end{equation}

where ``c'' is an arbitrary constant of the domain 
(\cite{goedel}, axiom III.1, p. 62).\\

This ``specialization rule'', as remarked before 
for the mathematical induction principle, has the logical deficiency of trying 
to express a metalogical concept within the object language. 
In fact, if we expand the antecedent of the sentence (\ref{eq14}), we 
find, among other conjunctive expressions, the expression to 
the right of the implication; in other words, we attain a trivial 
tautology. But the real purpose can only be described by a metalogical 
sentence which asserts that, if the formula $(\forall x) P(x)$ 
can be deduced then the formula P(x) can be deduced for any particular value of x (see formula (\ref{eq13})).\\

Combining the two sentences (\ref{eq12}) and (\ref{eq13}) we finally find:

\begin{equation}
(\forall x) Bew(P(x)) \leftrightarrow Bew( (\forall x) 
P(x) )
\label{eq15}
\end{equation}

Later it will be explained how a couple of sentences, postulated 
by G\"{o}del in the course of the proof of the undecidability theorem, 
and which predicate its own indeducibility, are in contradiction 
with the set of axioms used in that proof, when putting them 
together with the current theorem, that is, sentences (\ref{eq12}) or (\ref{eq15}).\\

The author would have preferred to propose rule (\ref{eq15}) as a metamathematical 
axiom, because it seems to have a more fundamental nature than 
the mathematical induction rule. However, in this paper, it is better to consider 
it as a theorem within the normal or current mathematics. It 
is also possible to propose another metamathematical axiom, similar 
to rule (\ref{eq15}), but using the existential quantifier instead of 
the universal quantifier.\\

\section{G\"{O}DEL'S UNDECIDABILITY THEOREM}

G\"{o}del's undecidability and incompleteness theorems impose the 
mathematicians the conclusion that the axiomatic methods have 
some intrinsic limitations that state, for example, that even 
the ordinary arithmetic cannot be fully axiomatized, or that 
most of the more significant fields of mathematics cannot be 
free of internal contradiction. If we can refute the limitative 
theorems, we could restore the bright alternatives proposed by 
Leibniz and Hilbert.\\

The first of G\"{o}del's limitative theorem, or undecidability 
theorem, has the number VI in the referenced author's original 
paper \cite{goedel}, since, to arrive to it, he shows a long development 
within the theory of ``primitive recursive functions''.\\

This theorem claims that, in system P (from Principia Mathematica 
augmented with Peano's axioms), there is always some sentence 
such that neither it nor its negation is deducible in the system.\\

The ``primitive recursive functions'' play a fundamental 
role in mathematics, since it is generally acknowledged that 
its use constitutes the formal equivalent of a ``finite 
effective method'' for computing or proving something; 
in other words, it means the same as what we are used to calling 
``algorithm''. The notion of mathematical truth has 
this character, so any human being is capable of reproducing 
a (mathematical) result.\\

If some theorem has been proved by using the tools of the ``primitive 
recursive functions'' theory, it is totally acceptable 
in principle; it cannot be said the same about ``non finitistic'' 
proofs which require the use of the ``infinite'' 
in order to prove something. \\

Preceding the proof of his theorem, G\"{o}del develops 46 relations 
and primitive recursive functions. For our current purposes we 
require only the last three definitions which are re-expressed 
in a little different notation; $\wedge$ between term represents 
logical conjunction, $\vee$ represents disjunction, $(\exists 
x)$ and $(\forall x)$ represent the existential and universal 
quantifiers respectively.\\

G\"{o}del's relation 46 expounds that a formula ``x'' 
is deducible if there exist a sequence ``y'' of formulas 
which represents the deduction of x:

\begin{equation}
deducible\_formula(x) \leftrightarrow
(\exists y) (y\ is\_a\_deduction\_of\_formula\ 
x)
\label{eq16}
\end{equation}

G\"{o}del's relation 45 is a definition of ``is\_a\_deduction\_of\_formula'' 
in terms of ``deduction'':\\

$y\ is\_a\_deduction\_of\_formula\ x \leftrightarrow $
\begin{equation}
deduction (y) \wedge x\ is\_last\_term\_of\ y
\label{eq17}
\end{equation}

His relation 44 is a definition of ``deduction'':\\

$deduction (y) \leftrightarrow (\forall x)( ( f\ is\_first\_term\_of\ 
y \wedge$

$l\ is\_last\_term\_of\ y \wedge$ 

$f \leq x \leq l) \rightarrow $

$is\_axiom (x) \vee (\exists p) (\exists q) (f \ensuremath{\leq} 
p,q \texttt{<} x \wedge$ 
\begin{equation}
x\ is\_an\_immediate\_inference\_of\ p,q) )
\label{eq18}
\end{equation}

We will not present here the definitions of ``is\_first\_term\_of'',\\ 
``is\_last\_term\_of'', ``is-axiom'' and 
``is\_an\_immediate\_inference\_of'', appealing to 
its intuitive interpretation or to the original G\"{o}del's paper \cite{goedel}. The three definitions mentioned explicitly above express that 
a formula is deducible if that formula is the last term of a 
deduction, where a deduction is a series of formulas in which 
each one of them is either an axiom or is obtained by an immediate 
inference of two preceding formulas.\\

Next, G\"{o}del sketches the proof of a theorem (theorem V in the 
original paper \cite{goedel}) which expresses that, for any n-ary relation 
R, it can be found a relational symbol ``r'' which 
can be deduced:

\begin{equation}
R(x1,..., xn) \rightarrow deducible\_formula (r(x1,..., xn))
\label{eq19.1}
\end{equation}

and also:

\begin{equation}
\sim R(x1,..., xn) \rightarrow deducible\_formula (\sim r(x1,..., 
xn))
\label{eq19.2}
\end{equation}

Later, G\"{o}del defines the relation:

\begin{equation}
Q(x, y) \leftrightarrow \sim x\ is\_a\_deduction\_of\_formula\ y 
\label{eq20}
\end{equation}

which is constructed by using only primitive recursive relations, 
so G\"{o}del concludes that Q(x, y) is also primitive recursive. 
This relation says that, for all ``x'' and for all 
``y'', x is not a deduction of formula y.\\

Consequently, by relation (\ref{eq20}), and theorem V, there must be a relational symbol ``q'' such that
\begin{equation}
\sim x\ is\_a\_deduction\_of\_formula\ y \rightarrow deducible\_formula\ 
(q(x,y))
\label{eq21}
\end{equation}
\begin{equation}
x\ is\_a\_deduction\_of\_formula\ y \rightarrow deducible\_formula 
(\sim q(x,y))
\label{eq22}
\end{equation}

Following this, G\"{o}del proposes the formula:

\begin{equation}
p = (\ensuremath{\forall} x)q(x,p)
\label{eq23}
\end{equation}

and calls: 

\begin{equation}
q(x,p) = r(x)
\label{eq24}
\end{equation}

which, after replacing in (\ref{eq23}), gives:

\begin{equation}
p = (\forall x)q(x,p) = (\forall x)r(x)
\label{eq25}
\end{equation}

This formula ``p'' declares about itself that it 
is not deducible. Take notice that the relation is not demonstrated 
and it is very strange to mathematics. It is in some sense similar 
to the equality relation
z = z where we postulate that z in the right hand side is an 
arbitrary function of itself, such as \makebox{f(x) = z+1}; after replacing in the equality, we get: z = f(z) 
= z+1. With this kind of operations we could prove anything.\\

Next, G\"{o}del substitutes ``y'' for ``p'' 
in (\ref{eq21}) and (\ref{eq22}), and applies the equivalences (\ref{eq24}) and (\ref{eq25}), 
obtaining:

\begin{equation}
\sim x\ is\_a\_deduction\_of\_formula\ (\forall x)r(x) \rightarrow 
deducible\_formula\ (r(x)) 
\label{eq26}
\end{equation}

\begin{equation}
x\ is\_a\_deduction\_of\_formula\ (\forall x)r(x) \rightarrow 
deducible\_formula\ (\sim\ r(x))
\label{eq27}
\end{equation}\\

Now G\"{o}del is prepared to prove the undecidability theorem (theorem 
VI in the original paper), which states that neither formula 
$(\forall x)r(x)$, nor formula $\sim(\forall x)r(x)$ 
are deducible in the system.\\

Thus, in the first part he proves in the following way that:\\

\textbf{{\LARGE 1. \ensuremath{\sim} deducible\_formula (}}\textbf{{\LARGE \ensuremath{\forall}}} \textbf{{\LARGE x)r(x)}}\\
For, if it were deducible, in other words if:\\
deducible\_formula (\ensuremath{\forall} x)r(x),\\
then, by (\ref{eq16}), there would be an ``n'' such that\\
n is\_a\_deduction\_of\_formula (\ensuremath{\forall} x)r(x).\\
But then, by (\ref{eq27}), \\
deducible\_formula (\ensuremath{\sim}r(n)),\\
Whereas, by the same:\\
deducible\_formula (\ensuremath{\forall} x)r(x)\\
it is also concluded that:\\
deducible\_formula (r(n)).\\
And the system would be inconsistent.\\

Next, in the second part, G\"{o}del proves that:\\

\textbf{{\LARGE 2. \ensuremath{\sim} deducible\_formula (\ensuremath{\sim}(}}\textbf{{\LARGE \ensuremath{\forall}}} 
\textbf{{\LARGE x)r(x))}}\\
It was just proved that:\\
\ensuremath{\sim} deducible\_formula (\ensuremath{\forall} x)r(x)\\
then, by (\ref{eq16}), \\
(\ensuremath{\forall} n) \ensuremath{\sim} n is\_a\_deduction\_of\_formula (\ensuremath{\forall} 
x)r(x).\\
From here it follows, by (\ref{eq26}), that\\
(\ensuremath{\forall} n) deducible\_formula (r(n)),\\
which, together with\\
deducible\_formula (\ensuremath{\sim}(\ensuremath{\forall} x)r(x) )\\
would be incompatible with the consistency of the system.\\
So, he concludes, (\ensuremath{\forall} x)r(x) is undecidable, finishing 
the proof of theorem VI.\\

In G\"{o}del's proofs in general, and in the present paper in particular, 
we should be aware that the arguments within the predicates are 
always G\"{o}del's numbers, which we have represented, for the 
sake of clarity, with the character series of the formulas they 
represent. The only operation allowed upon those arguments is 
the substitution of ``free variables'' when the unification 
is done in order to apply resolution.\\

\section{DISCOVERING THE CONTRADICTION}

Let us put together the set of rules, sentences or postulates 
used by G\"{o}del in the proof of theorem VI, not including yet 
the mathematical induction rule (\ref{eq12}). The names of the predicates will be 
abridged by employing the original G\"{o}del's terminology, where 
``deducible\_formula'' is represented by ``Bew'' 
and\\
``is\_a\_deduction\_of\_formula'' is represented by ``B''.

From (\ref{eq16}):
\begin{equation}
Bew(x) \leftrightarrow (\exists y) (y B x)
\label{eqn16}
\end{equation}

From (\ref{eq26}):
\begin{equation}
\sim x B (\forall x)r(x) \rightarrow Bew(r(x))
\label{eqn26}
\end{equation}

From (\ref{eq27}):
\begin{equation}
x B (\forall x)r(x) \rightarrow Bew(\sim r(x))
\label{eqn27}
\end{equation}

To this group of three clauses we have to add three sentences 
that G\"{o}del uses implicitly, without having proved or written 
previously. In part 1 of the proof of his theorem VI, he uses 
the two rules:

\begin{equation}
Bew (\forall x)r(x) \rightarrow Bew(r(y))
\label{eq28}
\end{equation}

\begin{equation}
\sim ( Bew(x) \wedge Bew(\sim x))
\label{eq29}
\end{equation}

Rule (\ref{eq28}) expresses that if a relation or property ``r(x)'' 
can be deduced for all values of variable x, then that relation 
can also be deduced for a particular value of ``y'' 
(G\"{o}del uses ``n'' as the value of ``y'').\\
Rule (\ref{eq29}) expresses the consistency of the system, in the sense 
that it cannot simultaneously be deduced a formula and its negation 
if the system is to be consistent.\\
In part 2 of his proof, G\"{o}del uses the following rule, also 
without having deduced or proved previously:

\begin{equation}
\ensuremath{\sim} ( (\ensuremath{\forall} y)Bew(r(y)) \ensuremath{\wedge} Bew(\ensuremath{\sim}(\ensuremath{\forall} 
x)Bew(r(x)) )
\label{eq30}
\end{equation}

which expresses that, if the system is going to preserve consistency, 
it is not possible to deduce ``r(y)'', for all values 
of ``y``, and also the negation of (\ensuremath{\forall} x)Bew(r(x)).\\
Now we re-express the premises in clausal form:\\

From (\ref{eqn16}16) are obtained

\begin{equation}
\sim Bew(x) \vee n B x
\label{eq31}
\end{equation}

\begin{equation}
Bew(x) \vee \sim y B x
\label{eq32}
\end{equation}

From (\ref{eqn26} 26): 
\begin{equation}
x\ B\ (\forall x)r(x) \vee Bew(r(x))
\label{eq33}
\end{equation}

From (\ref{eqn27}27): 
\begin{equation}
\sim x\ B\ (\forall x)r(x) \vee Bew(\sim r(y) )
\label{eq34}
\end{equation}

From (\ref{eq28}28): 
\begin{equation}
\sim Bew(\forall x)r(x) \vee Bew(r(y) )
\label{eq35}
\end{equation}

From (\ref{eq29}29): 
\begin{equation}
\sim Bew(x) \vee \sim Bew(\sim (x) )
\label{eq36}
\end{equation}

and, from (\ref{eq30}30): 
\begin{equation}
\sim Bew(r(n1)) \vee \sim Bew(\sim (\forall x)r(x) )
\label{eq37}
\end{equation}

Skolem constants n and n1 have been used to eliminate the existential quantifiers.\\

Let us call ``S'' to this set of clauses. By using 
the set S of clauses it is easy to reproduce G\"{o}del's proof 
by logical resolution:\\

To prove that (\ensuremath{\forall} x)r(x) is not deducible, we assume that it is deducible:

\begin{equation}
Bew(\forall x)r(x)
\label{eq38}
\end{equation}\\

resolving:\\

\begin{equation}
\ref{eq38} \wedge \ref{eq31}: n\ B\ (\forall x)r(x)
\label{eq39}
\end{equation}

\begin{equation}
\ref{eq39} \wedge \ref{eq34}:\ Bew(\sim r(y) )
\label{eq40}
\end{equation}

\begin{equation}
\ref{eq40} \wedge \ref{eq36}:\ \sim Bew( r(y) )
\label{eq41}
\end{equation}

\begin{equation}
\ref{eq38} \wedge \ref{eq35}:\ Bew( r(y) )
\label{eq42}
\end{equation}

\begin{equation}
\ref{eq41} \wedge \ref{eq42}:\ EMPTY\ CLAUSE
\nonumber
\end{equation}

As there is contradiction, it has been proved that:

\begin{equation}
\sim Bew(\forall x)r(x)
\label{eq43}
\end{equation}

From here, the proof can be followed with:

\begin{equation}
\ref{eq43} \wedge \ref{eq32}:\ 
\sim y\ B\ (\forall x)r(x)
\label{eq44}
\end{equation}

\begin{equation}
\ref{eq44} \wedge \ref{eq33}:\ Bew( r(x) )
\label{eq45}
\end{equation}

\begin{equation}
\ref{eq45} \wedge \ref{eq37}:\ \sim Bew(\sim (\forall x)r(x) 
)
\label{eq46}
\end{equation}

With this, G\"{o}del finishes the proof of theorem VI, for he wanted 
to come to sentences (\ref{eq43}) and (\ref{eq46}).

Now, sentence (\ref{eq45}) is interpreted by G\"{o}del as:

\begin{equation}
(\forall x)Bew(r(x))
\label{eq47}
\end{equation}

which declares that, for all x, the relation r(x) is deducible.\\

In the following, let us consider the theorem we had obtained 
from the rule of induction (sentence (\ref{eq12})):

\begin{equation}
(\forall x) Bew(r(x)) \rightarrow Bew (\forall x)r(x)
\label{eq48}
\end{equation}

Resolving between (\ref{eq47}) and (\ref{eq48}) we deduce:

\begin{equation}
Bew (\forall x)r(x)
\label{eq49}
\end{equation}

and the inconsistency is revealed, precisely between sentences 
(\ref{eq43}) and (\ref{eq49}).\\

It would be a meager favor to our formal systems if similar contradictions 
were to appear when we introduce the mathematical induction rule (\ref{eq48}). But 
we have to remark that the set S includes two special G\"{o}del's 
clauses, besides the rest which are perfectly acceptable and 
valid for any formal system. Those special clauses are (\ref{eq33}) and 
(\ref{eq34}), or similarly (\ref{eq26}) and (\ref{eq27}), which are arbitrarily postulated 
by G\"{o}del and which do not need to be present in our formal 
systems. If clauses (\ref{eq33}) and (\ref{eq34}) are suppressed then it is not 
anymore possible to prove the contradiction among the rest of 
clauses belonging to the set S, which now includes the rule of induction.\\

\section{CONCLUSION}

One of the corollaries of the undecidability theorem, the incompleteness 
theorem, establishes that the axiomatization, of any formal system 
which contains, at least, the elementary arithmetic, cannot be 
completed, unless it becomes inconsistent. The influence of this 
theorem in computer science lies in the fact that a computer 
program is directly expressible in, or translatable to, logic 
(for example PROLOG), which is a formal system. The answer, to 
the question of what would happen if, by accident, we had completed 
the axiomatization of arithmetic within a program, is that we 
could get absolutely whatever response, since that is what happens 
when logical contradiction is present in a formal system. So 
we would always have the doubt about the reliability of computation.\\

The dilemma of either stop using the mathematical induction rule 
or stop accepting G\"{o}del's theorem is established. It is suggested 
to use a very natural rule within the formal systems, which is 
the rule of induction, to avoid the effect of G\"{o}del's theorems.\\

\section{APPENDIX}

The following short note has the purpose of revealing the logical origin 
of the mathematical induction principle. This provides some useful insight and generalization for its use. 

\section{Definitions}

The mathematical induction principle is usually expressed in 
the following way:

Let P be any property of natural numbers.

Suppose that:\\

Q) 0 (zero) has the property P.

R) If any natural number x has property P, then its successor 
s(x) also has property P.

Then:

S) All natural numbers x have property P.\\

This ``mathematical induction principle'' can be 
expressed more concisely by using the logical connectives for 
conjunction ( \ensuremath{\wedge} ) and implication ( \ensuremath{\rightarrow} ) and 
the universal quantifier ( \ensuremath{\forall} ), as follows:

\begin{equation}
P(0)\ \wedge (\forall x)( P(x) \rightarrow P(s(x)) ) \rightarrow 
(\forall x) P(x)\ 
\label{a1}
\end{equation}

This formula is referred to by some authors as an ``axiom'', 
because it corresponds to the fifth Peano's axiom for the numbers (See G\"{o}del \cite{goedel}, p.61).

\section{The \textit{normal} use of the mathematical induction 
principle}

With the previous definitions, the form of the induction schema 
is: 
\makebox{ Q \ensuremath{\wedge} R \ensuremath{\rightarrow} S}\\

If the induction principle is to be applied to some problem, 
you have to ascertain that it is true Q (the base case) and, 
independently, it must be true R (the induction step). This means that you have the following three clauses:\\
1.\tab 
Q \ensuremath{\wedge} R \ensuremath{\rightarrow} S (induction principle)\\
2.\tab Q (base case)\\
3.\tab R (induction step)\\

Then, by using resolution between 1, 2 and 3, it is concluded 
``S''. This is the \textit{classic} or \textit{normal} use of 
induction.

\section{The logical origin of the mathematical induction principle}

If the mathematical induction principle were in fact an ``axiom'', 
it could not be demonstrated. But, in the following, it will 
be suggested an argument which seems to be a proof of it.\\

What I am going to do is to take only clauses 2 and 3, namely 
Q and R, or the base case and the induction step, and apply 
resolution between them. Obviously, we cannot use only the propositional 
letters but the full predicate notation. Let us assume that we 
apply the mathematical induction principle to the set of non-negative integers:

\begin{equation}
P(0)
\label{a2}
\end{equation}

\begin{equation}
(\forall x)( P(x) \rightarrow P(s(x)) )
\label{a3}
\end{equation}

Next, assign to x the value 0 in (\ref{a3}). It can be proved, by means 
of a truth table that the conjunction of clause (\ref{a2}) and the 
cited instantiation of (\ref{a3}) is logically equivalent to:

\begin{equation}
P(0) \wedge P(s(0))
\label{a4}
\end{equation}

With a new value for x in (\ref{a3}), now of s(0), clauses (\ref{a3}) and (\ref{a4}) 
are logically equivalent to:

\begin{equation}
P(0) \wedge P(s(0)) \wedge P(s(s(0)))
\label{a5}
\end{equation}

Expression (\ref{a5}) can be recognized that is an expansion of:

\begin{equation}
(\forall x) P(x)
\label{a6}
\end{equation}

where the universal quantifier ``for all x`` refers 
to the elements of the set: 

\begin{equation}
\{ 0,\ s(0),\ s(s(0)) \}
\label{a7}
\end{equation}

This is a growing set of elements. For, in each resolution step, 
a new element appears, and the cardinality of the set is increased 
by one. So, it is concluded that the left part of the mathematical 
induction principle, expressions (\ref{a2}) and (\ref{a3}) is logically equivalent to (\ensuremath{\forall}x)P(x), expression (\ref{a6}),
where the universal quantifier applies to a potentially infinite 
set of numbers. \\

This is the origin of the mathematical induction principle.

\section{What to think about this}

A ``potential infinite'' is, by definition, a ``growing 
but never completed set of elements''; and, according to 
Gauss, this is the only allowable infinite, because he didn't 
thought that the infinite could be closed.\\

In the previous section it has been shown that the left hand side of the mathematical induction 
principle implies the same result in the right 
hand side: (\ensuremath{\forall}x)P(x).\\

Consequently, we have found a tautology of the form: (\ensuremath{\forall}x)P(x) \ensuremath{\rightarrow} 
(\ensuremath{\forall}x)P(x)\\

But many current mathematicians do not accept this.\\

Because, when they conclude the right hand side, (\ensuremath{\forall}x)P(x), 
they are probably thinking in a closed infinite, as they are 
cantorian. \\

Moreover, the mathematicians feel ``entitled to stipulate 
what the domain is'', in the right hand side. \\

I see no reason to imply that the conclusion applies to some 
other set of elements, such as a set of characters, a set of 
lists, etc., different from the set of numbers with which I was working in the left hand side.\\

However, in some interesting mails received by me, this is precisely 
the objection given by a mathematician from the University 
of Illinois, to whom I will refer in the following by his initials 
``HD''. He concludes that the mathematical induction 
principle cannot be a theorem (or a tautology).\\

Let us analyze his example: Assume that the predicate ``P'' 
means ``even'' and that the function ``successor'' 
is defined as s(s(x)) (or ``x+2''); then, if we apply 
the induction principle, we obtain that the left hand side of 
the induction principle (namely Q \ensuremath{\wedge} R) is true; but, 
according to HD, the conclusion isn't true because not every 
number is even; or, in other words, the clause (\ensuremath{\forall}x)P(x) 
isn't true.\\

Let us analyze it in more detail. The individual resolutions between 
the base case and the induction step have shown us what the 
set of elements we are dealing with is. In fact, repeating the above 
argument, each instantiation of x, in each induction step (resolution), 
is increasing by one element the set of elements to which induction 
applies. Consequently, the interpretation of ``(\ensuremath{\forall}x)'' 
should be: ``all the elements for which either the base 
case applies (zero in the example) or the numbers obtained in 
each induction step''. We obtain the set: \{0, s(s(0)), 
s(s(s(s(0)))),...\}, which is the set of even numbers, to which 
the conclusion does apply. For this set it is true that (\ensuremath{\forall}x)P(x).\\

But, according to HD:\\

``You are right that the set \{0, s(s(0)), etc.\} is the 
set of even numbers. But the domain of my model is the set of 
all natural numbers. I am entitled to stipulate what the domain 
is. In that event, the universal quantifier says ``all 
natural numbers'' and the formula (\ensuremath{\forall}x)P(x) says 
``all natural numbers are even.'' I would add that 
it is well-known among mathematical logicians that the induction 
schema is falsifiable. HD''\\

What are the implications of this kind of reasoning? \\

First, it seems that the mathematicians jump to one conclusion 
that is not logically given by the premises. For, although in 
the previous example it is possible to continue indefinitely 
producing resolutions until the property for any given number 
is proved, the mathematicians conclude the same thing, and maybe 
faster, but for a different set of elements (a closed infinite) 
and, probably, because the god given law called mathematical 
induction principle says so.\\

Moreover, they are going a further step by supposing that, in 
the right hand side, they are entitled to stipulate what the 
domain is, with which the logical conclusion is not only extrapolated but lost. \\

With the use of the mathematical induction rule, which is a little different to what has already been said, and which has a metamathematical character, the author has been able to develop the refutation to the G\"{o}del incompleteness theorems, contained in the main body of this paper.\\

Such refutation consists in reexpressing, the left hand side 
of the mathematical induction rule (see Kleene \cite{kleene}), by using the \textit{classical} mathematical induction principle that has been reviewed in this appendix. It is obtained the following formula, which is logically equivalent to the induction rule:\\

(\ensuremath{\forall}x)Bew(P(x)) \ensuremath{\rightarrow} Bew((\ensuremath{\forall}x)P(x))\\

Where ``(\ensuremath{\forall}x)Bew(P(x))'' means that there 
are proofs for all: $\vdash$P(0), $\vdash$P(1), $\vdash$ P(2), etc.\\

``Bew'' is a predicate used by G\"{o}del with the meaning 
of ``is deducible''.\\

This formula enters in contradiction with the premises of G\"{o}del incompleteness and undecidability theorems, with the consequence that such theorems are invalid.\\

\end{document}